\documentclass[11pt,a4paper,subeqn,reqno]{article}

\usepackage{graphicx}
\usepackage{amssymb}
\usepackage{epstopdf}

\usepackage{latexsym}

\usepackage{amsmath,amsthm}

\usepackage[mathscr]{eucal}
\usepackage[all]{xy}
\usepackage{mathrsfs}
\usepackage{authblk}

\usepackage{amsfonts}
\usepackage{subfigure}
\usepackage{caption}
\usepackage[hyperpageref]{backref}
\usepackage[dvipsnames,usenames]{color}
\usepackage{color}
\usepackage{shuffle}
\usepackage{colortbl}
\usepackage{makeidx}
\usepackage[all]{xy}
\usepackage{shuffle}

\usepackage{datetime2}

\allowdisplaybreaks
\usepackage[text={152mm,225mm},centering]{geometry}

\usepackage{mathtools}
\usepackage{dsfont}
\usepackage{url}
\usepackage{enumerate}

\usepackage[hyperindex=true, colorlinks=false]{hyperref}

\newtheorem{theorem}{Theorem}

\newtheorem{lemma}{Lemma}[section]

\newtheorem{corollary}[lemma]{Corollary}

\theoremstyle{definition}

\newtheorem{example}[lemma]{Example}

\newtheorem{definition-lemma}[lemma]{Definition-Lemma}
\newtheorem{definition-theorem}[lemma]{Definition-Theorem}

\newtheorem{remark}[lemma]{Remark}

\newcommand{\Czero}{\mathbb{C}\{\zeta\}}
\newcommand{\Cone}{\mathbb{C}\{\zeta-1\}}

\newcommand{\Log}{\frac{\log(1-\zeta)}{2\pi i}}
\newcommand{\mI}{\mathcal{I}}

\newcommand{\Om}{\Omega}
\newcommand{\ii}{^{-1}}

\newcommand{\C}{\mathbb{C}}

\newcommand{\begla}{\begin{equation}}
\newcommand{\beglab}[1]{\begin{equation}	\label{#1}}
\newcommand{\edla}{\end{equation}}

\newcommand{\defeq}{\coloneqq}

\newcommand{\ens}{\enspace}
\newcommand{\ie}{{\it i.e.}\ }

\newcommand{\imp}{\ens \Rightarrow \ens}

\newcommand{\unit}{\frac1{1-\zeta}}

\newcommand{\datestampa}{{\small{File:\ens\hbox{\tt\jobname.tex}
\ens \today}}} 
\newcommand{\datestamp}{{\small{File:\ens\hbox{\tt\jobname.tex}
\ens \DTMnow}}} 

\begin{document}

\thispagestyle{empty}
\begin{center}
{\bf \Large
On the Hadamard inverse of a resurgent function with only one singularity}

\vspace{.7cm}

{Y.~Li\footnote{Yau Mathematical Sciences Center, Tsinghua University,
    Beijing},
  D.~Sauzin\footnote{CNRS -- Observatoire de Paris, PSL Research
    University, Paris, and Capital Normal University, Beijing}
  and S.~Sun\footnote{Capital Normal University, Beijing\\ \datestamp}
}

\bigskip 

{Preliminary version, October 2023
}

\end{center}
\smallskip

\begin{abstract}
  This note discusses the location of the singularities of the
  Hadamard inverse of an endlessly continuable function, in the case when
  the original function has only one singular singularity which is
  either a single pole or a simple singularity.
\end{abstract}
\bigskip


The Hadamard product of two power series is defined by termwise
multiplication of their coefficients (\cite{Ha98},\cite{Bo98}):
\begin{equation}
    \sum\limits_{n\geq 0} a_n \zeta^n \odot \sum\limits_{n\geq 0} b_n \zeta^n = \sum\limits_{n\geq 0} a_nb_n\zeta^n.
\end{equation}
With this product, the space $\Czero$ of  all holomorphic germs at the
origin becomes a commutative associative algebra, with unit $\delta :=
\frac{1}{1-\zeta}$.
In this note, we discuss a simplified ``Hadamard inverse
problem''. Namely, we consider the question: \emph{if $F \in \Czero$ has a
Hadamard inverse $G=F^{\odot-1}\in\Czero$ and if the analytic continuation of~$F$
has only one singular point, at which the singularity is either a pole
or a simple singularity, what are the
possible singular points of the analytic continuation of~$G$?}
\ \\

Here, as is usual in \'Ecalle Resurgence Theory (\cite{E81},
\cite{E85}, \cite{MS16}), we call ``simple singularity'' the superposition of a simple
pole and a logarithmic singularity with regular monodromy variation:
in the second case, our holomorphic germ~$F(\zeta)$ is assumed to have analytic continuation
along any path of~$\C$ that starts at~$0$ and avoids a certain
$\omega\in\C^*$, with
\[
F(\omega+\xi) = \frac{c}{\xi} + V(\xi)\frac{\log\xi}{2\pi i} + R(\xi),
\]
where $c\in\C$ and $V(\xi), R(\xi)\in\C\{\xi\}$.
A more general notion is that of ``simply ramified singularity'', in
which the term $\frac{c}{\xi}$ is replaced by an arbitrary polynomial
in $\frac{1}{\xi}$.
\ \\

Note that if the analytic continuation of~$F$ has more than one
singular points, then~$F^{\odot-1}$ may have infinitely many singular points.
An elementary example is given by 
\begin{equation}
  F(\zeta) \defeq 1+\frac{1}{1-\zeta} - \frac2{2-\zeta} \imp
  F^{\odot-1}(\zeta) = 1 + \sum\limits_{m\geq0} \frac{\zeta}{2^m-\zeta}.
\end{equation}
A much wilder example is given in \cite{BM92}:
For any $q,p\in\C^*$ such that $|q|<1$, $|p|=1$ and $p$ not a root of
unity,
$F(\zeta)\defeq \frac1{1-\zeta} - \frac{q}{1-p\zeta}$ has Hadamard
inverse $F^{\odot-1}(\zeta) = \sum\limits_{n\geq 0 } \frac{1}{1-qp^n} \zeta^n$,
which is shown to converge in the unit disk \emph{with the unit circle as
natural boundary.}
\ \\

On the other hand, we have proved in \cite{LSS20} that the Hadamard
product leaves invariant the space of \emph{resurgent functions}.
This space, pertaining to Resurgence Theory, is defined as follows:
Given a non-empty closed discrete subset~$\Om$ of~$\C$, one calls
``$\Om$-continuable germ'' any $f(\zeta)\in\Czero$ that admits
analytic continuation along any path of $\C\setminus\Om$ that starts in the disc of convergence
of~$f$;
one calls ``resurgent germ'' any $f(\zeta)\in\Czero$ that
is $\Om$-continuable for some~$\Om$.
\ \\

In this note we will focus on the case where~$F$ is assumed to have a
Hadamard inverse $G=F^{\odot-1}$ and to have only one singularity,
located at a point $\omega\in\C^*$.
\begin{itemize}
\item
In Section~\ref{sectionpole}, we assume that this singularity is a
pole and prove Theorem~\ref{propositionFpole}, to the effect that~$G$
is $\{0,\omega\ii\}$-continuable.
%
%
\item
  In Section~\ref{sectionFsimplesingularity}, we assume that $F$ has a
  simple singularity at~$\omega$, that $G$ is endlessly continuable,
  and that, at any $\omega'\in\C^*$, any branch of~$G$ has at worst a
  simple singularity.
  Then our conclusion is again that~$G$ is
  $\{0,\omega\ii\}$-continuable (Theorem
  \ref{theoremimplesingularity}).
\item
Section~\ref{sectionintegralformulaandproofs} gathers the proofs of
some auxiliary results stated in Section~\ref{sectionFsimplesingularity}.
  \end{itemize}
\medskip
  
\noindent Without loss of generality, we assume that the only singularity
 of~$F$ is located at~$1$.
  
\section{Hadamard inverse of germ with single pole}\label{sectionpole}

In this section, we assume $F$ is meromorphic in~$\C$, with a pole at
$1$ and no other poles. We will show that the Hadamard inverse
$G=F^{\odot-1}$ satisfies a linear non-homogeneous ordinary differential equation with
holomorphic coefficients, that is singular at worst at $0$
and~$1$. This will imply that the only possibly singular points of $G$ are $0$ and $1$. Let us see two examples first. \ \\
\begin{example}\label{examplesecondorder}
If $F=\frac{1}{(1-\zeta)^2}$, then we compute $F =  \sum\limits_{n\geq 0 } (n+1) \zeta ^n$ and $G= F^{\odot-1}\sum\limits_{n\geq0} \frac{1}{n+1} \zeta^n = -\frac{1}{\zeta} \log(1-\zeta)$. This means
  \begin{equation}
    \frac{1}{(1-\zeta)^2} \odot \left(-\frac{1}{\zeta} \log(1-\zeta) \right)  =\delta.
  \end{equation}
We notice that $G\in\Czero$ is resurgent function holomorphic at origin and it has singular points at $0$ and $1$. 
\end{example}

\begin{example}\label{examplethirdorder}
If $F=\frac{2}{(1-\zeta)^3}+\frac{1}{(1-\zeta)^2}$, then we compute $F =   \sum\limits_{n\geq 0 } (n+2)(n+1) \zeta ^n+\sum\limits_{n\geq0} (n+1)\zeta^n =\sum\limits_{n\geq0}(n+1)(n+3)\zeta^n$ and thus $G= \sum\limits_{n\geq0} \frac{1}{(n+1)(n+3)} \zeta^n =  \left(-\frac{1}{\zeta} \log(1-\zeta)\right) \odot \left(\frac{1}{\zeta^3}(-\log(1-\zeta)-\zeta-\frac{\zeta^2}{2}))\right)$. This means
  \begin{equation}
    \left(\frac{2}{(1-\zeta)^3}+\frac{1}{(1-\zeta)^2}\right) \odot \left(\frac{1}{\zeta} \log(1-\zeta)\right) \odot \left(\frac{1}{\zeta^3}\log(1-\zeta)+
    \frac{1}{\zeta^2}+\frac{1}{2\zeta}\right)  =\delta.
  \end{equation}
By a recent study about the resurgent property of Hadamard theorem (\cite{LSS20}), $G$ is a resurgent function with possibly singular points at $\{0,1\}$. 
\end{example}
%
%
%
\ \\

If $F$ has the form of $F_1(\zeta)+h(\zeta)$ where $h$ is an entire function and $F_1$ is resurgent which provides the singular part of $F$, then the following lemma shows that the singular behavior of $G=F^{\odot-1}$ is independent on the choice of $h$. \ \\
%
%
\begin{lemma}\label{lemmaentirenoeffect}
  Let $F,G$ be two germs holomorphic at origin with $F\odot G=\unit$. If $(F+h) \odot G_1 =\unit$ holds for an entire function $h$, then the difference of $G_1$ and $G$ is an entire function.
\end{lemma}

\begin{proof}
  Let \begin{equation}
    F(\zeta) = \sum\limits_{n\geq0} a_n \zeta^n, \quad     G(\zeta) = \sum\limits_{n\geq0} b_n \zeta^n, \quad
    h(\zeta) = \sum\limits_{n\geq0} c_n \zeta^n, \quad
    G_1(\zeta) = \sum\limits_{n\geq0} d_n \zeta^n,
  \end{equation}
  with $a_nb_n=1$, $(a_n+c_n)d_n=1$ for any $n\geq0$ and there exists $A,B,C$ s.t. for $n$ large enough and fixed $\epsilon>0$
  \begin{equation}
    |b_n| \leq AB^n, \quad |c_n| \leq C\epsilon^n.
  \end{equation}
  Then we compute
  \begin{equation}
    d_n=\frac{1}{a_n+c_n}=\frac{1}{\frac{1}{b_n}+c_n} = \frac{b_n}{1+c_nb_n} = b_n(1-\frac{c_nb_n}{1+c_nb_n}) = b_n-\frac{c_nb_n^2}{1+c_nb_n}.
  \end{equation}
  This induces that
  \begin{equation}
    G_1(\zeta) = G(\zeta) - \sum\limits_{n\geq0} \frac{c_nb_n^2}{1+c_nb_n} \zeta^n,
  \end{equation}
  where the last term in RHS is an entire function since for sufficient big $n$, we have an estimate
  \begin{equation}
    \left|\frac{c_nb_n^2}{1+c_nb_n}\right| \leq CA^2(\epsilon B^2)^n , \quad \epsilon \text{ small enough.}
  \end{equation}
\end{proof}


By using above lemma, we may assume 
\begin{equation}
    F(\zeta) = \sum_{j=1}^M \frac{a_j}{(1-\zeta)^j} \quad \text{with } a_M \neq 0
\end{equation} 
without an entire function and consider the singularity of $G=F^{\odot-1}$ in the latter discussion. The following lemma is useful.

\begin{lemma}\label{lemmaforinverseofapole}
 Let $f$ and $g$ be two holomorphic germs at origin. We have
  \begin{equation}
    \left( \partial_{\zeta}^sf \right)\odot g=\left(\sum\limits_{n\geq0}a_{n+s}\zeta^n\right) \odot\left(\partial_\zeta^s(\zeta^sg)\right).
  \end{equation}
  where $f(\zeta)=\sum\limits_{n\geq0}a_n\zeta^n$.
\end{lemma}

\begin{proof}
  Let $g(\zeta) = \sum\limits_{n\geq0} b_n\zeta^n$. We compute
  \begin{equation}\begin{split}
    \left( \partial_{\zeta}^sf \right)\odot g
    &=
    \left(\sum\limits_{n\geq s}(n+s)\cdots(n+1)a_{n}\zeta^{n-s}\right)\odot g=
    \left(\sum\limits_{n\geq0}(n+s)\cdots(n+1)a_{n+s}\zeta^n\right)\odot g
    \\
    &=
    \sum\limits_{n\geq0} (n+s)\cdots(n+1)a_{n+s}b_n\zeta^n
    =
    \left(\sum\limits_{n\geq0}a_{n+s}\zeta^n\right)
    \odot\left(\sum\limits_{n\geq0}(n+s)\cdots(n+1)b_n\zeta^n\right)
    \\
    &=
    \left(\sum\limits_{n\geq0}a_{n+s}\zeta^n\right) \odot\left(\partial_\zeta^s(\zeta^sg)\right).
  \end{split}\end{equation}
\end{proof}



We are ready to prove the first result we claimed in the introduction now. \ \\

\begin{theorem}\label{propositionFpole}
  If meromorphic function $F$ has only pole at $\omega\neq0$ and $G:=F^{\odot-1}$ is holomorphic at origin, then $G$ satisfies a differential equation which has singular points at most at $0$ and $\omega^{-1}$. Thus $G$ is a resurgent function with possibly singular points at $0$ and $\omega^{-1}$ in all sheets.
\end{theorem}

\begin{proof}
  By lemma \ref{lemmaentirenoeffect}, we directly assume  $F(\zeta)=\sum\limits_{j=1}^M\frac{a_j}{(\omega-\zeta)^j}=\sum\limits_{j=1}^M \frac{a_j}{(j-1)!}\partial_{\zeta}^{j-1}\frac{1}{\omega-\zeta}$. By Lemma \ref{lemmaforinverseofapole}, we have, for $G=\sum\limits_{n\geq0} b_n\zeta^n$,
  \begin{equation}\begin{split}
    \left( \partial_{\zeta}^s \frac{1}{\omega-\zeta} \right)\odot G
    &=
    \sum\limits_{n\geq0} (n+s)\cdots(n+1)\omega^{-n-s-1}b_n\zeta^n
    \\
    &=\partial_{\zeta}^s \left(\omega^{-s-1}\zeta^sG(\frac\zeta\omega)\right) = \left.\omega^{-1}\partial_X^{s} \left(X^sG(X)\right)\right|_{X=\frac\zeta\omega}.
  \end{split}\end{equation}
 A directly computation yields
  \begin{equation}
    \begin{split}
    F\odot G(\zeta)
      &=
      \left.\sum\limits_{j=1}^M \frac{a_j}{(j-1)!} \omega^{-1} \partial^{j-1}_X(X^{j-1}G(X))\right|_{X=\zeta/\omega}
      \\
      &=
      \left.\sum\limits_{j=1}^M \frac{a_j}{(j-1)!} \omega^{-1} \sum\limits_{s=0}^{j-1} \binom{j-1}{s} \frac{(j-1)!}{(j-1-s)!} X^{j-1-s}\partial_{X}^{j-1-s}G\right|_{X=\zeta/\omega}
    \end{split}
  \end{equation}
  Suppose $F \odot G = \delta$. It turns out that $G$ should satisfy the following differential equation:
  \begin{equation}
    \sum\limits_{j=0}^{M-1} c_j X^j \partial_X^j G = \frac{1}{1-\omega X}, \quad \text{with} \ c_{M-1}\neq0,
  \end{equation}
  which has only singular points at $0$ and $\omega^{-1}$. More precisely, if $M$ is $1$, the singular point of $G$ is $\omega^{-1}$; if $M\geq2$, the singular points of $G$ are $0$ and $\omega^{-1}$.
\end{proof}

\section{Hadamard inverse of a germ with single simple singularity}\label{sectionFsimplesingularity}

In this section, we are going to prove the following theorem.
\begin{theorem}\label{theoremimplesingularity}
   Assume $F,G$ are both holomorphic germs in $\Czero$ with $F\odot G =\delta$ and they are both resurgent. Assume $F,G$ have at most simple singularity at all points except at origin. If $F$ has only one singular point at $\omega$, then $G$ has possibly singular points at $0$ and $\omega^{-1}$ in all sheets.
\end{theorem}
\ \\
We always assume $\omega=1$ without lose of generality. A germ $\phi$ has a simple singularity at $1$ means that  
\begin{equation}
    F(\zeta) = \frac{A}{1-\zeta} + f_1(\zeta) \frac{\log(1-\zeta)}{2 
 \pi i} + f_2(\zeta)
\end{equation}
as $\zeta$ near $1$ where $ f_1,f_2 \in \Cone$. To prove the Theorem \ref{theoremimplesingularity}, it is easier to assume 
\begin{equation}\label{equationF}
    F(\zeta) = \frac A{1-\zeta} + f_1(\zeta) \frac{\log(1-\zeta)}{2 
 \pi i}  \quad  \text{ with }  A\neq 0
\end{equation}
where $f_1$ is an entire function in the light of Lemma \ref{lemmaentirenoeffect}. We will discuss the case as $A=0$ in Remark \ref{remarkA=0} latter. \ \\
\ \\
Let us see some lemmas at the beginning.
\begin{lemma}\label{lemmapoleodotpole}
    If $f_0(\zeta)=\sum\limits^M_{j=1}\frac{A_j}{(1-\zeta)^j}$, $g_0(\zeta)=\sum\limits^N_{j=1}\frac{B_j}{(\omega-\zeta)^j}$ with $A_M\neq0$, $B_N\neq0$, then
    \begin{equation}
        (f_0 \odot g_0)(\zeta)=\sum \limits^{M+N-1}_{j=0} \frac{C_j}{(\omega-\zeta)^j} \quad \text{with} \ C_{M+N-1}\neq0 .
    \end{equation}
\end{lemma}

\begin{lemma}\label{lemmapoleodotlog}
    If $f_0(\zeta)=\sum\limits^M_{j=1}\frac{A_j}{(1-\zeta)^j}$, $g\in \Czero$ resurgent and $g= g_1(\zeta)\frac{log(\omega-\zeta)}{2\pi i}+g_2(\zeta) \in \Czero$ with $g_1,g_2 \in \C\{\zeta-\omega\}$ in its principle sheet, then 
    \begin{equation}
        \lim\limits_{\zeta\rightarrow1}(\zeta -\omega)^M \left(f_0 \odot g \right) = 0. 
    \end{equation}
\end{lemma}

\begin{lemma}\label{lemmalogodotlog}
    If $f_1$ is an entire function, $g\in \Czero$ resurgent and $g= g_1(\zeta)\frac{log(1-\zeta)}{2\pi i}+g_2(\zeta) \in \Czero$ with $g_1,g_2 \in \mathbb{C}\{\zeta-\omega\}$ in its principle sheet, then
    \begin{equation}
    \begin{split}
        \lim\limits_{\zeta\rightarrow1}(\zeta -\omega) \left( \left( f_1(\zeta) \frac{log(1-\zeta)}{2 \pi i} \right) \odot g \right) = 0.
    \end{split}
    \end{equation}
\end{lemma}
\ \\
The proofs of these lemmas can be found in section \ref{sectionintegralformulaandproofs}. Due to a explicit formula in \cite{PM20a}, a directly corollary of the last lemma is as follows.

\begin{corollary}
    If $f_1$ is an entire function, $g\in \Czero$ resurgent and $g= g_1(\zeta)\frac{log(\omega-\zeta)}{2\pi i}+g_2(\zeta) \in \Czero$ with $g_1,g_2 \in \mathbb{C}\{\zeta-\omega\}$ in its principle sheet, then
    \begin{equation}\label{equationh1h2}
    \begin{split}
        \left( \left( f_1(\zeta) \frac{log(1-\zeta)}{2 \pi i} \right) \odot g(\zeta) \right) = h_1(\zeta) \frac{log(\omega-\zeta)}{2 \pi i} + h_2(\zeta).
    \end{split}
    \end{equation}
where 
\begin{equation}\label{equationMarco}
    h_1(\zeta) = -\frac1{2\pi i} \int_{\omega}^\zeta f_1(\frac{\zeta}{u}) g_{1}(u) \frac{du}{u} \in \mathbb{C}\{\zeta-\omega\}
\end{equation}
and $h_2 \in \Cone $.
\end{corollary}

\begin{corollary}
      Let $ F(\zeta) = \frac A{1-\zeta} + f_1(\zeta) \frac{\log(1-\zeta)}{2 
 \pi i} $. If $G\in \Czero$ is resurgent and $G$ has a simple singularity at $\omega$ in principle sheet, i.e.,
 \begin{equation*}
      G(\zeta) = \frac{B}{\omega-\zeta} + g_1(\zeta) \frac{\log(\omega-\zeta)}{2 \pi i} + g_2(\zeta)
 \end{equation*}
$g_1,g_2\in\mathbb{C}\{\zeta-\omega\}$,  then $F\odot G$ is singular at $\omega$.
\end{corollary}

\begin{proof}
    We prove it by contradiction. By the order estimate in previous lemmas, one may find that the highest order of the polar part of $F\odot G$ at $\omega$ is provided by $ \frac{AB}{\omega-\zeta}$. Thus either $A=0$ or $B=0$ (or maybe both). Let $B=0$. Using formula \eqref{equationh1h2} and \eqref{equationMarco}, if $F\odot G$ has no monodromy at $\omega$, then 
    \begin{equation}
        Ag_1(\zeta) - \frac1{2\pi i} \int_{\omega}^\zeta f_1(\frac{\zeta}{u}) g_{1}(u) \frac{du}{u} =0.
    \end{equation}
    One may compute the derivatives of the left hand side of above equation and take the limit as $\zeta$ goes to $\omega$. Thus, either $g_1(\zeta)=0$ or $f_1(\zeta)=0$ for any $A$. Another way to see it is to use the method about unique solution to Volterra integral equation \cite{VV}, \cite{Br04}. 
\end{proof}
\ \\
Now let us fix the form of resurgent function $G$ in the principle sheet:
\begin{equation}\label{equationG}
           G(\zeta) = g_0(\zeta) + g_1(\zeta) \frac{\log(1-\zeta)}{2 \pi i} + g_2(\zeta),
\end{equation}
where
\begin{equation}
g_0 = \frac{B}{1-\zeta}, \text{ with } B\neq 0 \quad g_1,g_2 \in\Cone.
\end{equation}
\begin{remark}
The reason why we choose a Hadamard inverse in above form is as
follows. Roughly speaking, if both of $F,G$ have simply ramified
singularities at~$1$, \ie 
\begin{equation}
    \begin{split}
        F(\zeta) &= f_0(\zeta) + f_1(\zeta) \frac{\log(1-\zeta)}{2  \pi i} \text{ with } f_0 = \sum\limits_{j=1}^M{ \frac{A_j}{(1-\zeta)^j}}, \ A_M\neq 0 
 \\
      G(\zeta) &= g_0(\zeta) + g, \text{ with } g_0 = \sum\limits_{j=1}^N{ \frac{B_j}{(1-\zeta)^j}}, \ B_N\neq 0 \text{ and } g=g_1(\zeta) \frac{\log(1-\zeta)}{2 \pi i} + g_2(\zeta)
    \end{split}
\end{equation}
as $\zeta$ near $1$ with $f_1,g_1,g_2\in\Cone$, then the orders of the polar parts of $f_0 \odot g_0$, $f_0 \odot g$ and $(f_1\log(1-\zeta)) \odot g$ at $1$ are exactly $M+N-1$, less than $M$ and less than $1$ correspondingly. Moreover, $F\odot g_2$ is holomorphic at $1$ due to Hadamard theorem. These facts imply that, as $M=1$ in Theorem \ref{theoremimplesingularity}, $N=1$ in our ansatz is a very reasonable choice. 
\end{remark}
\ \\
\begin{remark}\label{remarkA=0}
    The reason for $A\neq0$ in the assumption of equation \eqref{equationF} is as follows. Let $A=0$. By Lemma \ref{lemmapoleodotlog} and Lemma \ref{lemmalogodotlog}, if $G$ has a simple singularity at $1$, then 
    \begin{equation}
        \lim\limits_{\zeta \rightarrow 1} (\zeta-1) (F\odot G) =0,
    \end{equation}
    which is a contradiction.
\end{remark}
\ \\
Now let $F,G$ be resurgent functions holomorphic at origin. Moreover, 
\begin{equation}
\begin{split}
     F&=\frac{A}{1-\zeta}+f_1(\zeta)\frac{log(1-\zeta)}{2\pi i},
     \\
      G&=\frac{B}{1-\zeta}+g_1(\zeta)\frac{log(1-\zeta)}{2\pi i}+g_2(\zeta)
    \end{split}
\end{equation}
$A,B\neq0$, $f_1$ entire and $g_1,g_2 \in \Cone$ as assumed in formulas \eqref{equationF} and \eqref{equationG}. The rest of this section is devoted to proving Theorem \ref{theoremimplesingularity} under this setting. By using formula \eqref{equationh1h2}, one finds
\begin{equation}
F\odot G=\frac{AB}{1-\zeta}+A\left(g_1\frac{log(1-\zeta)}{2\pi i}+g_2\right)+Bf_1(\zeta)\frac{log(1-\zeta)}{2\pi i}+h_1(\zeta) \frac{log(1-\zeta)}{2 \pi i} + h_2(\zeta) 
\end{equation}
and further
\begin{equation}\label{equationFodotGconditions}
     F\odot G =\frac{1}{1-\zeta} 
    \Longrightarrow
    \begin{cases}
        A \cdot B = 1 \quad\Longrightarrow\quad B=A^{-1}
        \\
        Ag_1 + Bf_1 + h_1 =0
        \\
        Ag_2 = -h_2
     \end{cases}.
\end{equation}

To prove the Theorem \ref{theoremimplesingularity} is sufficient to prove that $g_1$ has only possibly singular point at $0$ and $g_2$ has possibly singular points at $0$ and $1$. Indeed, by the classical theory of Volterra integral equation (see for instance \cite{VV} or \cite{Br04}), the second equation in right hand side of \eqref{equationFodotGconditions} determines a unique solution $g_1$ with only possibly singular points at $0$. 
The story for $g_2$ is much more complicated because of  less information. We prove that $g_2$ has only possibly singular points at $0$ and $1$ by contradiction. Notice that $g_2$ satisfies
\begin{equation}\label{equationconditionforg2}
    f_1(g)\frac{log(1-\zeta)}{2\pi i}\odot\left(g_1\frac{log(1-\zeta)}{2\pi i}+g_2\right)=h_1\frac{log(1-\zeta)}{2\pi i}-Ag_2.
\end{equation}
with known $f_1,g_1,h_1$. We are going to prove that if $g_2$ has simple (even simply ramified) singularity at $\omega\neq 0,1$ in a sheet corresponding to a path $\gamma$ which start at origin, then above formula no longer holds. In fact, one may discuss the singular behavior of the both sides of formula \eqref{equationconditionforg2} at $\omega$ and prove the difference of them. The following lemmas focus on the order of the polar parts of different Hadamard products as $\zeta$ goes to $\omega$ along the path $\gamma$.

\begin{lemma}\label{lemmag2polaratomega}
Let $f_1$ be an entire function and $F(\zeta)=f_1(\zeta)\frac{log(1-\zeta)}{2\pi i}$. Let $g$ have simply ramified singularity at $\omega$ associated with a path $\gamma$: $(cont_\gamma g)(\zeta)=\sum\limits_{j=1}^M\frac{A_j}{(\omega-\zeta)^j}+\phi_1(\zeta)\frac{log(\omega-\zeta)}{2\pi i}+\phi_2(\zeta)$ with $A_M\neq0$ and $\phi_1,\phi_2 \in\mathbb{C}\{\zeta-\omega\}$ as $\zeta$ near $\omega$. Then 
\begin{equation}
    \lim\limits_{\zeta\rightarrow\omega}(\zeta-\omega)^M cont_\gamma(F\odot g)(\zeta)=0.
\end{equation}
\end{lemma}
\ \\
This lemma implies that there is no polar part in $g_2$ at $\omega$ associated to a path $\gamma$. Indeed, if we assume $g_2$ has a pole of order $M$ at $\omega$ , then the left hand side of formula \eqref{equationconditionforg2} has order less than $M$ but the right hand side has order exactly $M$, which is a contradiction. Let us define further the singular behavior of $g_2$ at $\omega$ associated to a path $\gamma$ is 
\begin{equation}
    (cont_\gamma g_2)(\zeta)=\phi_1(\zeta)\frac{log(\omega-\zeta)}{2\pi i}+\phi_2(\zeta) \quad \text{ with }\phi_1,\phi_2 \in\mathbb{C}\{\zeta-\omega\}.
\end{equation}

\begin{lemma}\label{lemmag2logatomega}
Let $f_1$ be an entire function and $F(\zeta)=f_1(\zeta)\frac{log(1-\zeta)}{2\pi i}$. Let the resurgent function $g$ have singularity at $\omega$ in the branch associated with $\gamma:(cont_\gamma g)(\zeta)=\phi_1(\zeta)\frac{log(\omega-\zeta)}{2\pi i}+\phi_2(\zeta)$ with $\phi_1,\phi_2 \in\mathbb{C}\{\zeta-\omega\}$. Then 
\begin{equation}
    \lim\limits_{\zeta\rightarrow\omega}(\zeta-\omega)cont_{\gamma}(F\odot g)(\zeta)=0.
\end{equation}.    
\end{lemma}
\ \\
Using again the monodromy formula in \cite{PM20a}, this lemma allows us to write the monodromy of the left hand side of \eqref{equationconditionforg2} to be 
\begin{equation}
    -\frac1{2\pi i} \int^{\zeta}_\omega f_1\left( \frac{\zeta}{u}\right) \phi_1(u) \frac{du}{u},
\end{equation}
which should coincide to the monodromy of the right hand side $\phi_1(\zeta)$. The theory of Volterra integral equation says the unique solution to  
\begin{equation}
 -\frac1{2\pi i} \int^{\zeta}_\omega f_1\left( \frac{\zeta}{u}\right) \phi_1(u) \frac{du}{u} +A\phi_1(\zeta) =0
\end{equation}
should be $\phi_1=0$. This means that $g_2$ has no singular behavior at $\omega$, which yields the desired result as in Theorem \ref{theoremimplesingularity}.
\ \\

\section{Appendix}\label{sectionintegralformulaandproofs}

\subsection{Integral formula of Hadamard product}        
        
\begin{center}
\begin{figure}[h]
\centering
\includegraphics[width=2in,height=2in]{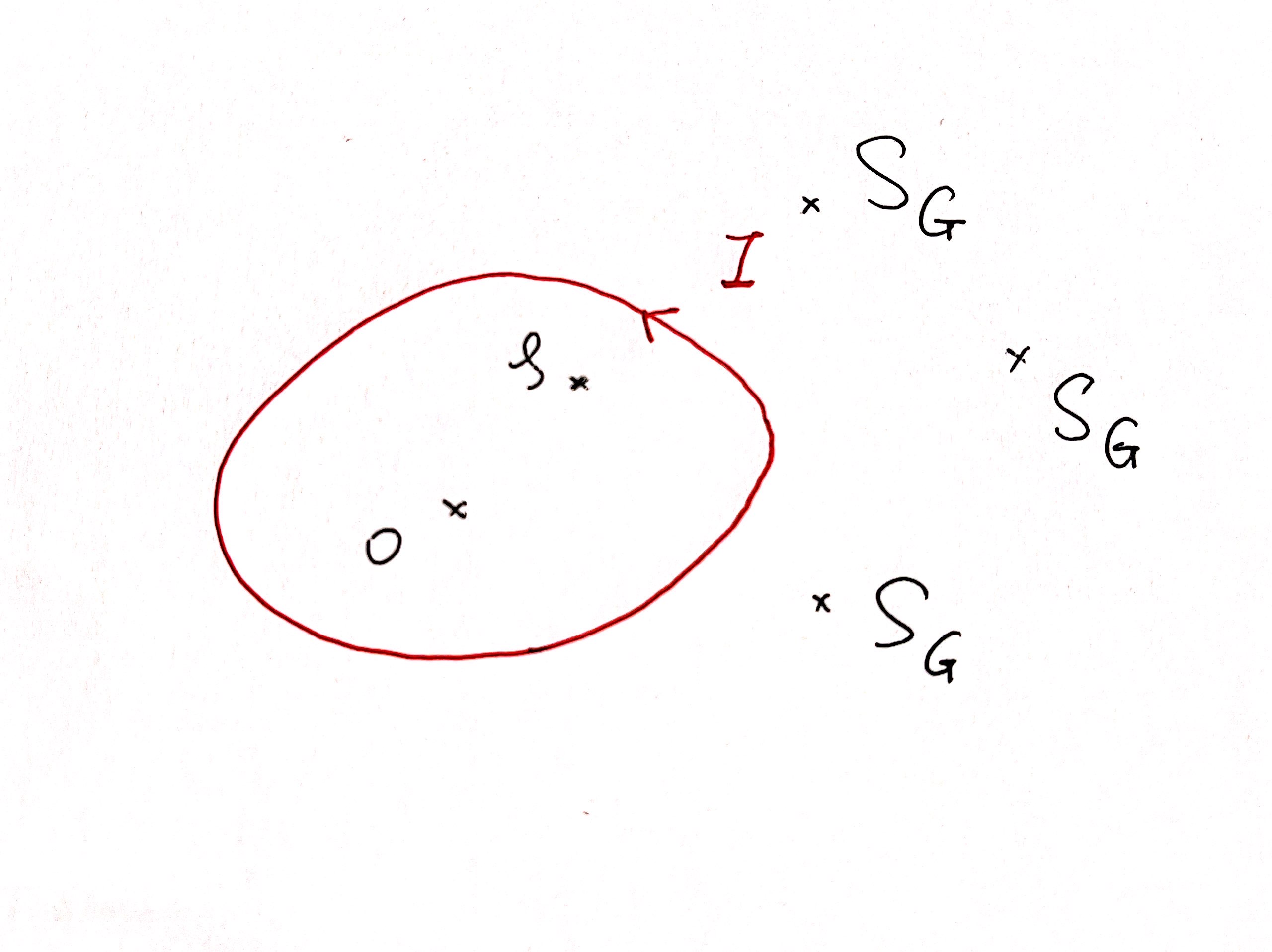}
\caption{ $\mathcal{I}$ is a small circle around both $0$ and $\zeta$ anticlockwise as $\zeta$ small enough.}
\label{curveI}
\end{figure}
\end{center}

\begin{center}
\begin{figure}[h]
\centering
\includegraphics[width=5.5in,height=2in]{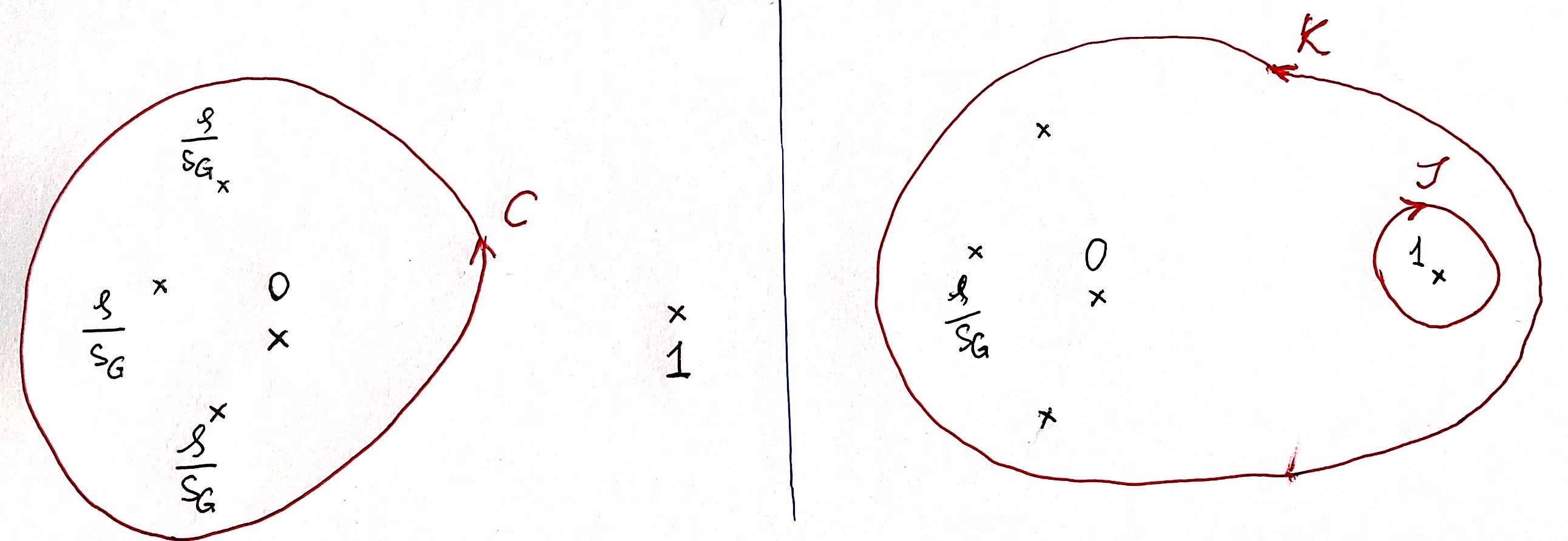}
\caption{LHS: $\mathcal{C}$ is a cirle with radius less than 1 which contains all points induced by the singular points of $G$. If $F$ has no branch at $1$, then $\mathcal{C}$ can be divided by $\mathcal{K}$ and $\mathcal{J}$ in RHS.}
\label{curveCKJ}
\end{figure}
\end{center}

In this section we will recall the integral formula of Hadamard product in the special cases required for this note. One may see classical articles \cite{Bo98}, \cite{Ha98} or recently \cite{LSS20} for a general version. Let $F \in \Czero$ and it has a single simply remified singularity at $1$, i.e.,
\begin{equation}
    F=\sum_{j=1}^M \frac{a_j}{(1-\zeta)^j} + f_1(\zeta) \Log  + f_2(\zeta)
\end{equation}
with entire functions $f_1,f_2$. Let $G$ be a resurgent function holomorphic at origin and it admits analytic continuation along any path which avoids $S_G\in\mathbb{C}$. The integral formula of Hadamard product for $\zeta$ small enough is
\begin{equation}\label{equationintegralonI}
  F \odot G (\zeta)  = \frac 1{2\pi i} \oint _{\mathcal{I}} F(\frac{\zeta}{z}) G (z) \frac{dz}{z} 
\end{equation}
where $\mathcal{I}$ is a small circle around both $0$ and $\zeta$ anticlockwise (see Figure \ref{curveI}).
On the other hand, by the commutative property of Hadamard product, one may exchange the positions of $F$ and $G$ 
\begin{equation}\label{equationintegralonC}
  F \odot G (\zeta)  = \frac 1{2\pi i} \oint _{\mathcal{C}} F(z) G (\frac{\zeta}{z}) \frac{dz}{z} 
\end{equation}
where $\mathcal{C}$ is a cirle with radius less than 1 which contains all points induced by the singular points of $G$. Moreover, if $F$ is meromorphic, i.e. $f_1=0$, then we may divide the integral along $\mathcal{C}$ by two parts
\begin{equation}\label{equationintegralonKJ}
    F \odot G (\zeta)  = \frac 1{2\pi i} \oint _{\mathcal{K}} F(z) G (\frac{\zeta}{z}) \frac{dz}{z} + \frac 1{2\pi i} \oint _{\mathcal{J}} F(z) G (\frac{\zeta}{z}) \frac{dz}{z}
\end{equation}
where $\mathcal{J}$ is a small curve around $1$ clockwise and $\mathcal{K}$ is an anticlockwise curve with large enough radius. See Figure \ref{curveCKJ}.

It is remarkable that the first integral (along $\mathcal{K}$) is an entire function in $\zeta$: as $\zeta$ goes along $\gamma$, the ``moving singular points'' of this integral run along $\frac{\gamma}{g_i}$ for $g_i\in S_G$, but they never pinch the homotopy class of $\mathcal{K}$. Thus, the singularity of $F\odot g$ is only provided by the second integral.

\subsection{Proofs}

\begin{proof}[Proof of Lemma \ref{lemmapoleodotpole}]
    We only prove in the case of $\omega=1$ for simplicity. The proof is similar for the other points on the principle sheet. By using formula \eqref{equationintegralonI} and residue formula, a directly computation yeilds
    \begin{equation}
    \begin{split}
      (f_0 \odot g_0)(\zeta)
    =& \sum \limits^M_{j=1} \sum \limits^N_{k=1} \frac{1}{2\pi i} \oint_{\mathcal{C}}\frac{A_j}{(1-\frac{\zeta}{z})^j}\frac{B_k}{(1-z)^k}\frac{dz}{z}
    \\
    =& \sum \limits^M_{j=1} \sum \limits^N_{k=1} \frac{A_j B_k}{2\pi i}\oint_{\mathcal{C}}\frac{z^{j-1}}{(z-\zeta)^j}\frac{1}{(1-z)^k}dz
    \\
    =& \sum \limits^M_{j=1} \sum \limits^N_{k=1} \frac{1}{(j-1)!}A_j B_k\left(\frac{d}{d\zeta}\right)^{j-1} \left(\zeta^{j-1}\frac{1}{(1-\zeta)^k} \right) .
    \end{split}
    \end{equation}
    Thus the order of the pole of $(f_0 \odot g_0)(\zeta)$ at $1$ is $M+N-1$, which is provided by $j=M$ and $k=N$.
\end{proof}

\begin{proof}[Proof of Lemma \ref{lemmapoleodotlog}] We would like to use two different ways to prove it.     We only prove in the case of $\omega=1$ for simplicity. The proof is similar for the other points on the principle sheet. \ \\
\textbf{Way 1:} Using formula \eqref{equationintegralonI} and residue formula, we compute 
 \begin{equation}
    \begin{split}
         &f_0\odot\left(g_1(\zeta)\frac{log(1-\zeta)}{2\pi i}+g_2(\zeta)\right) 
        \\
        =& \frac{1}{2\pi i}\oint_{\mathcal{C}}\sum\limits^M_{j=1}\frac{A_j z^{j-1}}{(z-\zeta)^j} 
        \left(g_1(z)\frac{log(1-z)}{2\pi i} + g_2(z)\right) dz
        \\
        =& \sum\limits^M_{j=1} A_j \frac{1}{(j-1)!} \left(\frac{d}{d\zeta}\right)^{j-1} 
        \left(\zeta^{j-1}\left(g_1(\zeta)\frac{log(1-\zeta)}{2\pi i}+g_2(\zeta)\right)\right).
    \end{split}
\end{equation}
Then the desired result holds since the highest order of the pole at $\zeta=1$ in above summation is $M-1$.\ \\
\ \\
\textbf{Way 2:} Let us use formula \eqref{equationintegralonKJ} since $f_0$ has no branch at $1$. As mentioned before, the first term in the right hand side of equation \eqref{equationintegralonKJ} is an entire function. Thus we only need to prove 
  \begin{equation}
    \lim\limits_{\zeta\rightarrow1} (\zeta-1)^M \oint_{\mathcal{J}} f_0(z) \left(g_1(\frac{\zeta}{z})\frac{1}{2\pi i} \log(1-\frac{\zeta}{z}) +g_2(\frac{\zeta}{z}) \right) \frac{dz}{z}=0.
  \end{equation}
The integral can be parametrized as follows. Let $\mathcal{J}=1+re^{i\theta}$ and $r=|\zeta-1|\sin\varphi$ (see Figure \ref{parametrizeJ}), then we have
\begin{equation}\begin{split}
  &\oint_{\mathcal{J}} f_0(z) \left(g_1(\frac{\zeta}{z})\frac{1}{2\pi i} \log(1-\frac{\zeta}{z}) +g_2(\frac{\zeta}{z}) \right) \frac{dz}{z}
  \\
  =&
  \frac{1}{2\pi}\int_0^{2\pi} f_0(1+re^{i\theta})  \left(g_1(\frac{\zeta}{1+re^{i\theta}})\log(1-\frac{\zeta}{1+re^{i\theta}})+ g_2(\frac{\zeta}{1+re^{i\theta}})\right) \frac{re^{i\theta}}{1+re^{i\theta}} d\theta.
  \end{split}
\end{equation}
Recall that $f_0$ has a pole at $1$ with order at most $M$ and $g_1,g_2$ holomorphic at $1$. We end our proof because of the following estimate:
\begin{equation}
  \left| (\zeta-1)^M f_0(1+re^{i\theta})g_1(\frac{\zeta}{1+re^{i\theta}})\log(1-\frac{\zeta}{1+re^{i\theta}}) \frac{re^{i\theta}}{1+re^{i\theta}}\right| \rightarrow 0, \quad\text{as } \zeta\rightarrow 1.
\end{equation}

\begin{center}
\begin{figure}[h]
\centering
\includegraphics[width=3in,height=1.5in]{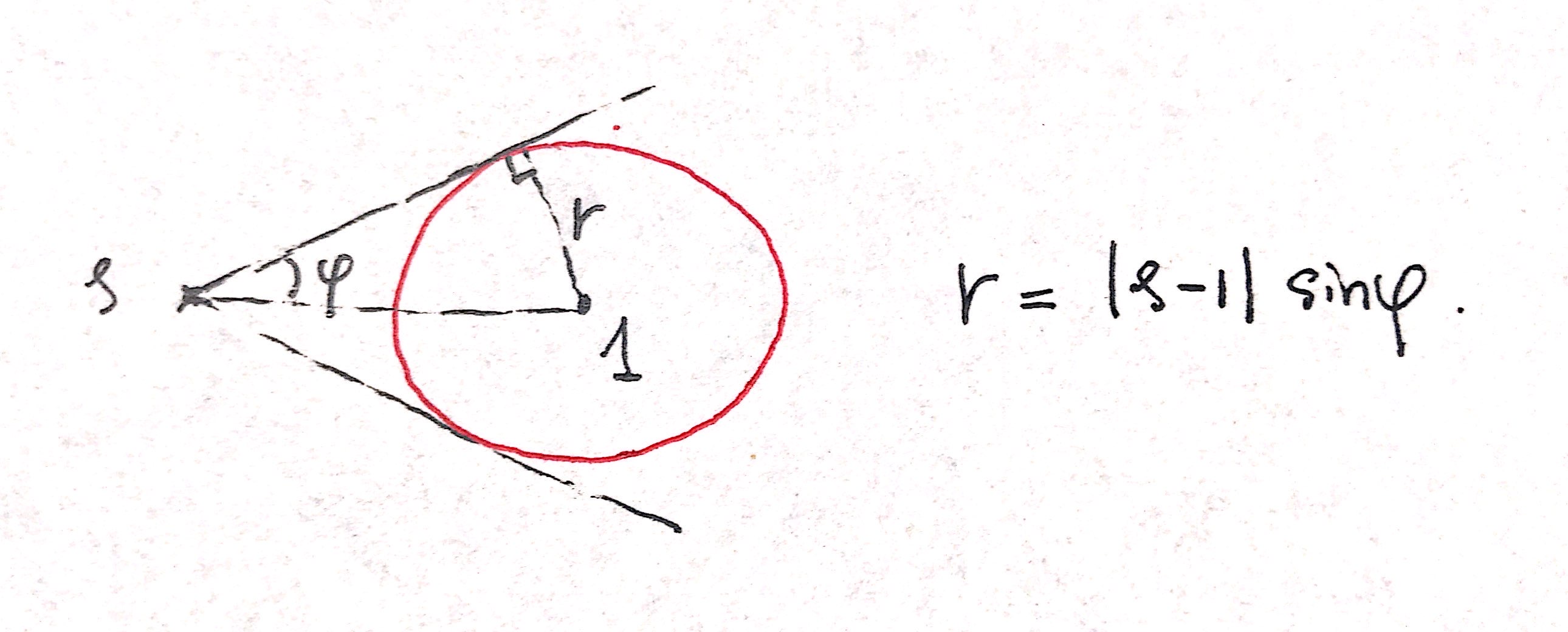}
\caption{$r$ and $|1-\zeta|$ perpetually maintain an unchanging proportion with fixed argument $\varphi$.}
\label{parametrizeJ}
\end{figure}
\end{center}
\end{proof}

\begin{proof}[Proof of Lemma \ref{lemmalogodotlog}]
    We only prove in the case of $\omega=1$ for simplicity. The proof is similar for the other points on the principle sheet. By using formula \eqref{equationintegralonI}, it is sufficient to prove

\begin{equation}
    \lim\limits_{\zeta\rightarrow1}(\zeta -1) \oint_{\mathcal{I}} 
     H(\zeta,z)dz = 0
\end{equation}
with 
\begin{equation}
    H(\zeta,z) :=\frac1z f_1(\frac{\zeta}{z}) \frac{\log(1-\frac{\zeta}{z})}{2\pi i}\left(g_1(z)\frac{log(1-z)}{2\pi i}+g_2(z)\right).
\end{equation}
\ \\
As shown in \cite{LSS20}, there exists a non-autonomous vector field s.t. its solution is a curve $\mathcal{I}^\prime$ homotopy to $\mathcal{I}$ as $\zeta$ goes to $1$ along $\gamma$. See the Figure \ref{logodotloganalyticcontinuation}.
\begin{center}
\begin{figure}[h]
\centering
\includegraphics[width=4in,height=1.2in]{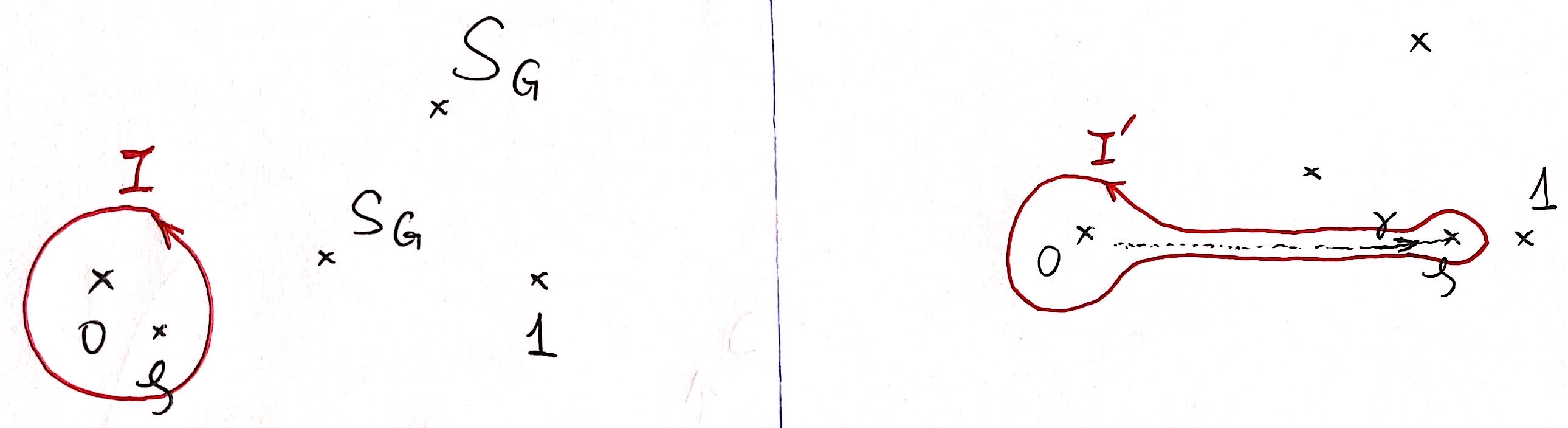}
\caption{}
\label{logodotloganalyticcontinuation}
\end{figure}
\end{center}
\ \\
Let's decompose the curve $\mathcal{I}^\prime=\mI_1+\mI_2+\mI_3+\mI_4$ as in Figure \ref{logodotlogdecomposition}. For simplicity, we assume $\zeta<1$. 
\begin{center}
\begin{figure}[h]
\centering
\includegraphics[width=5in,height=1.5in]{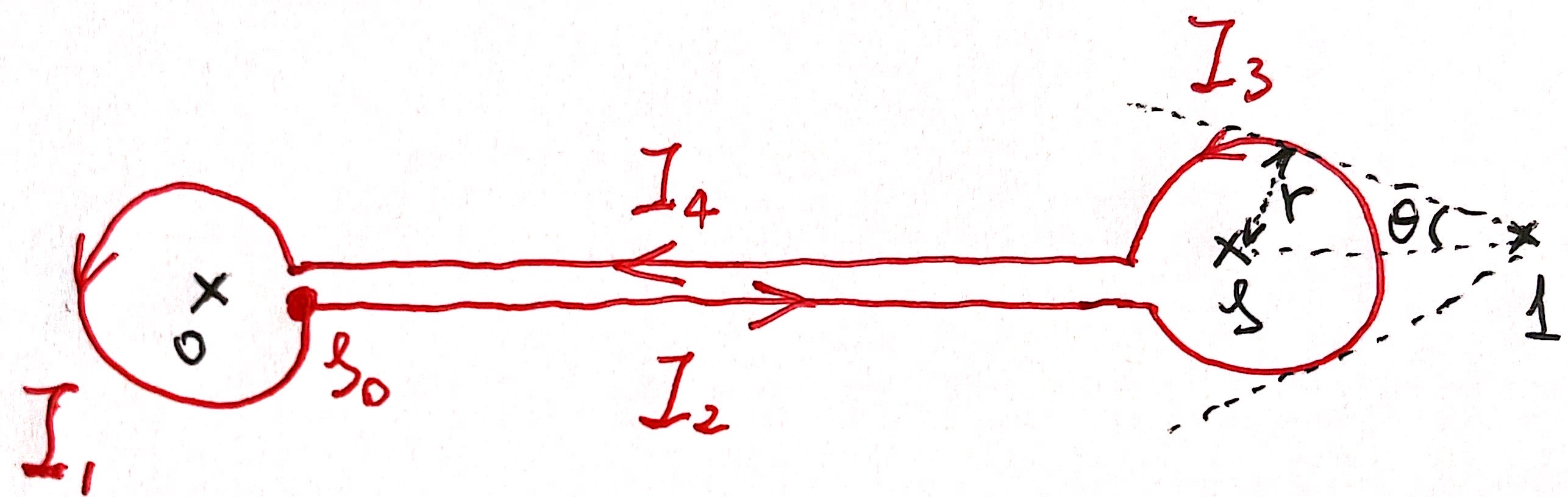}
\caption{For a fixed small $\zeta_0$, we denote the $\mathcal{I}_1$ by the curve goes around $0$ anticlockwise starting and ending at $\zeta_0$. Let $\mI_2$ and $\mI_4$ be two curves with same ending points at $\zeta_0$ and $\zeta-r$ but different directions. Let $\mI_3$ be a small curve goes around $\zeta$ anticlockwise with radius $r$. As $\zeta$ goes to $1$, $r = (1-\zeta) \sin\theta$ goes to $0$ with a fixed $\theta$. }
\label{logodotlogdecomposition}
\end{figure}
\end{center}
One can prove that 
\begin{equation}
    \int_{\mI_1} H(\zeta,z)dz
\end{equation}
is bounded since both $\mI_1$ and $H$ are bounded. Moreover,

\begin{equation}
\begin{split}
     \int_{\mI_2+\mI_4} H(\zeta,z)dz 
     =&
     \int_{\zeta_0}^{\zeta-r} H(\zeta,z)dz + \int_{\zeta-r}^{\zeta_0} H(\zeta,z)+  f_1(\frac{\zeta}{z})g_1(z)\frac{\log(1-z)}{2\pi i}\frac{dz}{z}
     \\
     =&
     -\frac{1}{2 \pi i} \int_{\zeta_0}^{\zeta-r} f_1(\frac{\zeta}{z})g_1(z)log(1-z)\frac{dz}{z}.
\end{split} 
\end{equation}
This is integrable as $\zeta \rightarrow1$ since $f_1$ is integrable at $1$. Finally, parametrizing $\mI_3$ in $ z=\zeta+re^{i \alpha},\alpha \in [-\pi,\pi], r=sin\theta|1-\zeta|$ with fixed $\theta$, 
%
one proves that $\lim\limits_{\zeta \rightarrow1} (1-\zeta)\int_{\mI_3} H dz = 0$ is implied by
\begin{equation}
    \left|log\left( 1-\frac{\zeta}{\zeta+re^{i\alpha}}\right)log(1-\zeta-re^{i\alpha})r\right|\rightarrow 0
\end{equation}
as $ \zeta\rightarrow1$. Thus,
\begin{equation}
    \lim\limits_{\zeta \rightarrow1}(1-\zeta) \int_{\mI^\prime} H dz =  
    \lim\limits_{\zeta \rightarrow1} (1-\zeta) \int_{\mI_1+\mI_2+\mI_3+\mI_4} H dz = 0.
\end{equation}
\end{proof}

The methods in the proofs for Lemma \ref{lemmag2polaratomega} and Lemma \ref{lemmag2logatomega} are similar to the previous. The key point to prove them is to find the integral curve homotopy to original circle on the logarithmic Riemann surface at $1$ instead of complex plane. We may omit them here and leave them as further study.

\end{document}